\documentstyle[12pt]{article}
\renewcommand{\baselinestretch}{1.3}

\makeatletter
\newcommand{\single}{\let\CS=\@currsize\renewcommand{\baselinestretch}{1}\tiny\CS}
\newcommand{\singles}{\let\CS=\@currsize\renewcommand{\baselinestretch}{1.3}\tiny\CS}
\newcommand{\oneanda}{\let\CS=\@currsize\renewcommand{\baselinestretch}{1.2}\tiny\CS}
\newcommand{\doubles}{\let\CS=\@currsize\renewcommand{\baselinestretch}{1.5}\tiny\CS}
\newcommand{\tree}{\let\CS=\@currsize\renewcommand{\baselinestretch}{1.5}\tiny\CS}
\newcommand{\four}{\let\CS=\@currsize\renewcommand{\baselinestretch}{2}\tiny\CS}

\newcommand{\ncom}{\newcommand}
\ncom{\bq}{\begin{equation}}
\ncom{\eq}{\end{equation}}
\ncom{\beqn}{\begin{eqnarray*}}
\ncom{\eeqn}{\end{eqnarray*}}
\ncom{\beq}{\begin{eqnarray}}
\ncom{\eeq}{\end{eqnarray}}
\ncom{\been}{\begin{enumerate}}
\ncom{\eeen}{\end{enumerate}}
\ncom{\nno}{\nonumber}
\ncom{\hs}{\mbox{\hspace{.25cm}}}
\ncom{\rar}{\rightarrow}
\ncom{\lrar}{\longrightarrow}
\ncom{\Rar}{\Rightarrow}
\ncom{\noin}{\noindent}
\newtheorem{thm}{Theorem}[section]
\newtheorem{lemma}[thm]{Lemma}
\newtheorem{cor}[thm]{Corollary}
\newtheorem{pro}[thm]{Proposition}
\newtheorem{example}[thm]{Example}
\newtheorem{remark}[thm]{Remark}
\newtheorem{defn}[thm]{Definition}

\ncom{\bd}{\begin{defn}}
\ncom{\ed}{\end{defn}}
\ncom{\bt}{\begin{thm}}
\ncom{\et}{\end{thm}}
\ncom{\bl}{\begin{lemma}}
\ncom{\el}{\end{lemma}}
\ncom{\bco}{\begin{cor}}
\ncom{\eco}{\end{cor}}
\ncom{\bp}{\begin{pro}}
\ncom{\ep}{\end{pro}}
\ncom{\bex}{\begin{example}}
\ncom{\eex}{\end{example}}
\ncom{\brm}{\begin{remark}}
\ncom{\erm}{\end{remark}}
\ncom{\comx}{I\!\!\!\!C}
\ncom{\zee}{$Z\!\!\!\!Z$}
\ncom{\ze}{Z\!\!\!\!Z}
\ncom{\Q}{$I\!\!\!\!Q$}
\ncom{\p}{I\!\!P}
\ncom{\al}{\alpha}
\ncom{\be}{\beta}
\ncom{\f}{\frac}
\ncom{\ga}{\gamma}
\ncom{\bib}{\bibitem}
\ncom{\pf}{{\bf Proof: }}
\ncom{\sta}{\stackrel}
\ncom{\cA}{{\cal A}}
\ncom{\cG}{{\cal G}}
\ncom{\cI}{{\cal I}}
\ncom{\cO}{{\cal O}}
\ncom{\cV}{{\cal V}}
\ncom{\cW}{{\cal W}}
\ncom{\cK}{{\cal K}}
\ncom{\cE}{{\cal E}}
\ncom{\cL}{{\cal L}}
\ncom{\cZ}{{\cal Z}}
\ncom{\cP}{{\cal P}}
\ncom{\cN}{{\cal N}}
\ncom{\s}{(\!\!\times}
\title{Projective normality of abelian varieties}
\author{Jaya N.Iyer\\FB6, Mathematik, \\ Universitat GH Essen,
\\45117, Essen, Germany\\(Email: bms051@sp2.power.uni-essen.de) }
\begin{document}

\maketitle

Mathematics Classification Number: 14C20, 14K05, 14K25, 14N05.

\section{Introduction}
Let $L$ be an ample line bundle of type $\delta=(d_1,d_2,...,d_g)$ on
an abelian variety $A$ of dimension $g$. Consider the associated
rational map $\phi_L:A\lrar \p H^0(L)$. 
Suppose $L=M^n$, for some ample line bundle $M$ on $A$. Then Koizumi
and Ohbuchi have shown that $L$ gives a projectively normal embedding if $n\geq 3$ and
when $n=2$, no point of $K(L)$ is a base point for $M$, (see [BL],
7.3.1).
Consider the case when $g=2$, and $L$ is an ample line bundle of type
$(1,d)$ on $A$, i.e. $L\neq M^n$ for any ample line bundle $M$ on
$A$, $n>1$. Then
it has been shown by Lazarsfeld (see [L]) that whenever $\phi_L$ is
birational onto its image and $d\geq 7$ odd and $d\geq 14$ and even,
then $\phi_L$ gives a projectively normal embedding.
We showed that if the $\it{Neron\,Severi\,group}$ of $A$, $NS(A)$, is $\ze$,
generated by $L$, and $d\geq 7$, then $\phi_L$ gives a projectively
normal embedding, (see [I]).

In this article we show
\bt
Suppose $L$ is an ample line bundle on a $g$-dimensional simple abelian
variety $A$. If $h^0(L)>2^g.g!$ then $L$ gives a projectively normal
embedding, for all $g\geq 1$.
\et

Since projective normality is an open condition, our theorem is, therefore
true for a generic pair $(A,L)$, as above.

We outline the proof of 1.1.
 
We firstly show that for abelian varieties, it is enough to show the surjectivity of the
homomorphism 
 $$Sym^2H^0(L)\sta{\rho_2}{\lrar} H^0(L^2)$$ to give a projectively
 normal embedding, (see 2.3). 
 
 To show $2$-normality and hence projective normality of $L$, we
 firstly consider a finite isogeny $A\lrar B=A/H$, where $H$ is a
 maximal isotropic subgroup of the fixed group of $L$, $K(L)$. Then
 $L$ descends down to a principal polarization $M$ on $B$.
 We then show that the surjectivity of the map $\rho_2$ is equivalent to showing
  the dual subgroup $H'$ of $H$, in $B(\simeq Pic^0(B))$ generates
 the linear system of $M^2$( and hence its translates also) i.e. the images of points of $H'$, under
 the morphism $B\sta{\phi_{t^*_\sigma M^2}}{\lrar} \cK(B)_\sigma\subset |t^*_\sigma M^2|$,
 $b\mapsto t^*_b\theta+ t^*_{-b+\sigma}\theta$, have their linear span as
 $|t^*_\sigma M^2|$, for all $\sigma\in H'$.
( Here $\theta$ is the unique divisor in $|M|$ and $\cK(B)$ ( respectively 
$\cK(B)_\sigma$) denotes
 the Kummer variety of $B$ in $|M^2|$ ( respectively in $|t^*_\sigma M^2|$). 

We show ( see 3.2)
\bp
Let $\cL$ be an ample line bundle on a simple abelian variesty $Z$ and
consider the associated rational map $Z\sta{\phi_{\cL}}{\lrar}
\p H^0(\cL)$.
Then any finite subgroup $G$ of $Z$, of order strictly greater than 
$h^0(\cL).g!$, generates the linear system $\p H^0(\cL)$. More precisely,
the points $\phi_{\cL}(g)$ where $g$ runs over all elements of $G$ not
in the base locus of $\cL$ span $\p H^0(\cL)$.
\ep

We then apply above proposition to $\cL= t^*_\sigma M^2$, to obtain bounds
as asserted for a polarized abelian variety $(A,L)$, in 1.1

 $\it{Acnowledgements}$: We thank A.Hirschowitz and M.Waldschmidt for a useful
conversation. We are grateful to Marc Hindry for a discussion which
helped us in 3.1 and the referee for suggestions. We also thank Institut de Mathematiques, Univ.Paris-6, for their hospitality where this work was done and the French Ministry of National
Education, Research and Technology and DFG 'Arithmetic und Geometrie' Essen, for their support.

$\bf{Notations:}$
Let $\cL$ be an ample line bundle on an abelian variety Z, of dimension $g$.

The $\it{fixed\, group}$ of $\cL$ is $K(\cL)=\{a\in A: \cL\simeq
t_a^*\cL\},\,t_a:A\lrar A,x\mapsto a+x$.

The $\it{theta\, group}$ of $\cL$ is
$\cG(\cL)=\{(a,\phi):\cL\sta{\phi}{\simeq}t_a^*\cL\}$.

The $\it{Weil\, form}\, e^\cL:K(\cL)\times K(\cL)\lrar \comx^*$, is the
commutator map $(x,y)\mapsto x'y'x'^{-1}y'^{-1}$, for any lifts
$x',y'\in \cG(\cL)$, of $x,y\in K(\cL)$. 

$h^0(\cL)= dim H^0(Z,\cL)$

If $G$ is a finite subgroup of $Z$, then $Card(G)= order (G)$.

\section{`r-normality'\, of \, L}

Consider an abelian variety $A$ of dimension $g$ and an ample line
bundle $L$ on $A$.

Consider the multiplication maps
$$H^0(L)^{\otimes r}\sta{\rho_r}{\lrar} H^0(L^r),\, for\, r\geq 2.$$

\bd
 $L$ is said to be $r$-normal, if $\rho_r$ is
surjective.
\ed
\bd 
$L$ is normally generated, if $L$ is $r$-normal for all $r\geq2$.
\ed

The main result of this section is the following.

\bp
Suppose $L$ is an ample line bundle on an abelian variety $A$. If $L$
is $2$-normal, then  $L$ is $r$-normal, for all $r\geq 2$. In
particular, $L$ is normally generated.
\ep

Firstly, we will see
\bp
Suppose $L$ and $M$ are ample line bundles on an abelian variety
$A$. 

1) The multiplication map 
 $$\sum_{\al\in U}H^0(L\otimes \al)\otimes H^0(M\otimes
 \al^{-1})\lrar H^0(L\otimes M)$$ is surjective, for any open subset
 $U$ of $Pic^0(A)$.

2) If the multiplication map $H^0(L)\otimes H^0(M)\lrar H^0(L\otimes M)$ is surjective, then
   the maps
$$(a)\,H^0(L)\otimes H^0(M\otimes \al)\lrar H^0(L\otimes M\otimes \al)$$
and $$(b)\, H^0(L\otimes \al^{-1})\otimes H^0(M\otimes \al)\lrar
H^0(L\otimes M)$$ are also surjective, for $\al$ in some open subset
$U$ of $Pic^0(A)$.
\ep
\pf 1) See [BL], 7.3.3.

2) Denote $\hat{A}=Pic^0(A)$. Consider the projections $p_A:A\times
   \hat{A}\lrar A$ and $p_{\hat{A}}:A\times \hat{A}\lrar \hat{A}$ and
   the sheaves
   $\cE_0=p_{\hat{A}*}(p_A^*L),\,\cE_1=p_{\hat{A}*}(p_A^*L\otimes
   \cP ^{-1})$, $\cE_2=p_{\hat{A}*}(p_A^*M\otimes
   \cP) $, where $\cP$ is the Poincare bundle on $A\times \hat{A}$.
   Since the fibres $\cE_0(\al)= H^0(L), \cE_1(\al)=H^0(L\otimes \al)$
   and $\cE_2(\al)=H^0(M\otimes \al^{-1})$ have constant dimension, for
   all $\al\in \hat{A}$, $\cE_0,\cE_1$ and $\cE_2$ are vector bundles
   on $\hat{A}$.

Consider the natural maps
$$\cE_0\otimes \cE_2 \sta{\rho_{02}}{\lrar }p_{\hat{A}*}(p_A^*(L\otimes
  M)\otimes \cP)$$
and
$$ \cE_1\otimes \cE_2\sta{\rho_{12}}{\lrar} H^0(L\otimes M)\otimes
\cO_{\hat{A}}$$

Since the map $$\rho_{02}(0)=\rho_{12}(0): H^0(L)\otimes H^0(M)\lrar
H^0(L\otimes M)$$ by assumption, is surjective, by semi-continuity,
$\rho_{02}(\al)$ and $\rho_{12}(\al)$ are surjective, for $\al$ in
some open subset $U$ of $\hat{A}$. $\Box$

$\bf{Proof \,of\, 2.3:}$
We prove by induction on $r$.
Suppose the map $\rho_r:H^0(L)^{\otimes r}\lrar H^0(L^r)$ is
surjective.

Consider the map 
$$H^0(L)^{\otimes r+1}\sta{Id\otimes \rho_{r}}{\lrar} H^0(L)\otimes
H^0(L^{r})\sta{\rho_{1,r}}{\lrar }H^0(L^{r+1}).$$

To see the surjectivity of the map $\rho_{r+1}=\rho_{1,r}\circ
(Id\otimes \rho_r)$, we need to
show that the map $\rho_{1,r}$ is surjective.

By 2.4 1),$$ H^0(L).H^0(L^r)=\sum_{\al\in U}H^0(L).H^0(L\otimes
\al^{-1}).H^0(L^{r-1}\otimes \al).$$ 
 Since $L$ is $2$-normal, by 2.4 2) (a),
$$H^0(L).H^0(L\otimes \al^{-1})= H^0(L^2\otimes \al^{-1})$$
which implies ( using 2.4 1)) that
$$H^0(L).H^0(L^r)=\sum_{\al\in U} H^0(L^2\otimes
\al^{-1}).H^0(L^{r-1}\otimes \al)$$
$$= H^0(L^{r+1}).$$
$\Box$ 

\section{`2-normality' of $L$}
Consider the multiplication map $$H^0(L)\otimes
H^0(L)\sta{\rho_2}{\lrar} H^0(L^2).$$
This map factors via
$$H^0(\p H^0(L),\cO(2))=Sym^2H^0(L)\sta{\rho_2}{\lrar}H^0(L^2)$$
and $Ker\rho_2=I_2=$ the vector space of quadrics containing
$\phi_L(A)$ in $\p H^0(L)$.

We will use the following.

\bp Let $D$ be an ample divisor on a $g$-dimensional simple abelian variety
$Z$. Suppose $G$ is a finite subgroup of $Z$, of order $k$ ($=
Card(G)$, in the sequel), and
contained in $D$. Then $k\leq D^g$. 
( Here $D^g$ denotes the self intersection number of $D$, which by
Riemann-Roch is $h^0(\cO (D)).g!)$.
\ep
\pf Let $G\subset Z$ be a subgroup of order $k$ with $G\subset D$.
This implies that $G\subset Y=\cap_{\sigma\in G}D+\sigma$, where
$D+\sigma=\{x+\sigma: x\in D\}$, and now $Y$ is invariant under the
group $G$. If $s=dim Y=0$, then our proof ends here. Otherwise, let
$Y=Y_1\cup ...\cup Y_r$, where $Y_j$ are the irreducible components,
with $s_j=dim Y_j$. Let $Y_1$ be an irreducible component of $Y$, with 
$dimY_1=dimY$. Choose $g-s$ translates $D+h_j$, $h_j\in Z$, of $D$ which
 intersect properly
with $Y$, i.e. $Y\cap_{j=1}^{g-s} D+h_j$ is a finite set of points. Then 
$Y_1\cap_{j=1}^{g-s} D+h_j \subset  Y\cap_{j=1}^{g-s} D+h_j$ and hence
 $degY_1\leq deg Y\leq D^g$. ( Here $deg Y$ is
the intersection number of $Y$ with the class of $D^{g-s}$, in
$H^*(Z)$).

Since $Y$ is $G$- invariant, $\cup_{\sigma\in G}Y_1+\sigma \subset Y$.
 Consider the subgroup $G_{Y_1}=\{g\in G:Y_1+g=Y_1\}$ of $G$. Since 
$\sum_{\sigma\in \frac{G}{G_{Y_1}}}deg(Y_1+\sigma)\leq deg Y$ and $deg Y_1=
deg Y_1+\sigma $ we get 
the inequalities $Card(\frac{G}{G_{Y_1}}).deg Y_1\leq deg Y\leq D^g$, 
i.e. $ Card(G)\leq D^g.\frac{Card(G_{Y_1})}{deg Y_1}$.
To complete our proof, it suffices to show that  $Card(G_{Y_1})\leq
deg Y_1$. Now $G_{Y_1}\subset Stab(Y_1)=\{a\in
Z: Y_1+a=Y_1\}$. Observe that $Stab(Y_1)=\cap_{y\in Y_1}Y_1-y$. Now for a 
$y_0\in Y_1$,
$Stab(Y_1)=(Y_1-y_0)\cap_{y\in Y_1,y\neq y_0} Y_1-y\subset (Y_1-y_0)\cap_{g\in G,y\in Y_1,y\neq y_0}D+g-y$. Hence it is now clear that $deg
Stab(Y_1)\leq deg Y_1$. Since $Z$ is simple, $Stab(Y_1)$ is
zero-dimensional and we get $Card(G_{Y_1})\leq deg Stab(Y_1)\leq deg
Y_1$.
 $\Box$
 
We will require the geometric interpretation of 3.1, namely,
\bp
Let $\cL$ be an ample line bundle on a simple abelian variety $Z$ and
consider the associated rational map $Z\sta{\phi_{\cL}}{\lrar}
\p H^0(\cL)$.
Then any finite subgroup $G$ of $Z$, of order strictly greater than 
$h^0(\cL).g!$, generates the linear system $\p H^0(\cL)$. More precisely,
the points $\phi_{\cL}(g)$ where $g$ runs over all elements of $G$ not
in the base locus of $\cL$ span $\p H^0(\cL)$.
\ep

$\bf{Proof\, of\, 1.1:}$

Consider a polarized simple abelian variety $(A,L)$, of type
$(d_1,...,d_g)$.
Let $H\subset K(L)$ be a maximal isotropic subgroup, for the Weil form
$e^L$.
Consider the isogeny $A\sta{\pi}{\lrar} B=\frac{A}{H}$. Then $L$
descends down to a principal polarization $M$ on $B$. We may assume
$M$ is symmetric , i.e, $M\simeq i^*M$, $i(b)=-b,b\in B$. By Projection
formula and using the fact that $\pi_*\cO_A= \oplus_{\chi\in \hat{H}}
L_{\chi}$, where $L_\chi$ denotes the degree $0$ line bundle on $B$,
corresponding to the character $\chi$ on $H$, we get
$$H^0(L)=\oplus_{\chi\in \hat{H}}H^0(M\otimes L_{\chi})= \oplus_{\sigma\in
  H'}H^0(t^*_\sigma M)$$
via the isomorphism $B\sta{\psi_M}{\lrar} Pic^0(B)$, $b\mapsto
t^*_bM\otimes M^{-1}$, $H'=\psi^{-1}_M(\hat{H})$.

Similarly, $$H^0(L^2)=\oplus_{\sigma\in H'}H^0(t^*_\sigma M^2).$$
 Consider the multiplication map $$Sym^2H^0(L)\sta{\rho_2}{\lrar} H^0(L^2).$$

Then $Sym^2H^0(L)=\sum_{\sigma,\tau\in
  H'}H^0(t^*_{\tau}M).H^0(t^*_{\tau^{-1}\sigma}M)$ and we can write the
  map $\rho_2$ as $\oplus_{\sigma\in H'}\rho_\sigma$, where
$$\sum_{\tau\in
  H'}H^0(t^*_\tau M).H^0(t^*_{-\tau+\sigma}M)\sta{\rho_\sigma}{\lrar}
  H^0(t^*_\sigma M^2)...(I).$$
To show surjectivity of $\rho_2$, it is enough to show surjectivity of
$\rho_\sigma$, for each $\sigma\in H'$.

Now it is well known that the morphism associated to the line bundle
$M^2$, embeds the Kummer variety of $B$, $\cK(B)$, in the linear system
$|M^2|$ and is given as 
$$ B\sta{\phi_{M^2}}{\lrar}\cK(B)\subset |M^2|, b\mapsto t^*_b\theta
+t^*_{-b}\theta$$
where $\theta$ is the unique symmetric divisor in $|M|$.

Under translation by $\sigma\in B$, we have the corresponding morphism 
$$ B\sta{\phi_{t^*_\sigma M^2}}{\lrar}\cK(B)_\sigma \subset |t^*_\sigma M^2|, b\mapsto t^*_b\theta
+t^*_{-b+\sigma}\theta.$$
( Here $\cK(B)_\sigma= Image(B)$ in $|t^*_\sigma M^2|$). 
So showing surjectivity of $\rho_\sigma$ in (I) is equivalent to
showing that the image of $H'$ in $\cK(B)$, generates the linear
system $|M^2|$ and hence the translates $|t^*_\sigma M^2|$. Since the pair $(A,L)$ is a simple
polarized abelian variety, with $h^0(L)= Card(H') > 2^g.g!=
h^0(M^2).g!$, by 3.2 and (I), each $\rho_\sigma$ is surjective.
Hence, by 2.3, our proof is now complete.

\brm
Notice that if $g=1$, any line bundle of degree strictly greater than
$2$ on an elliptic curve, gives projectively normal embedding. Hence
the bound is sharp.

If $g=2$ and $L(1,d)$ is an ample line bundle on an abelian surface $A$,
with $h^0(L)=7, 8$, then by [I] and [L], we know $L$ gives
projectively normal embedding (generically). By the method of proof of
1.1, we cannot expect to obtain a sharp bound since if L were to be of
type $(2,4)$, then it does not give a projectively normal embedding.

Consider the situation when $g\geq 3$.  
It is clear, by Kunneth, that if $A=A_1\times A_2\times...\times A_r$ and
$L=p^*_1L_1 \otimes...\otimes p^*_rL_r$, where $(A_j,L_j)$ are
polarized abelian varieties and $L_j$
give projectively normal embedding then $L$, of type $\delta$, also
gives a projectively
normal embedding. Hence for a generic pair $(A,L)$ of type $\delta$,
$L$ gives a projectively normal embedding.
This will show that there exists line bundles L with $h^0(L)\leq
2^g.g!$ on $A$, such that $A$ is projectively normal in $\p H^0(L)$.

But we can hope to improve the bound in 3.1, for $g>1$, by taking into consideration
the structure of the finite subgroup $H'$, which essentially
distinguishes
the type $\delta$ of $L$, for instance, types $(1,8)$ and $(2,4)$,
when $g=2$. 
Here $H'\simeq \frac{\ze}{8\ze}$ in the former case and in the latter
case either $\simeq \frac{\ze}{2\ze}\times \frac{\ze}{4\ze}$ or
$\simeq (\frac{\ze}{2\ze})^3$( satisfying certain condition with
respect to the Weil form $e^{\cL^2}$).

Also, in 3.1, if the ample divisor $D\subset Z$ is moreover a
symmetric divisor, i.e. $i(D)=D$, then
the subset $Y=\cap_{g\in G}D+g$ is in fact invariant for the action of
$G\times <i>$. So one may get better bounds in some cases and hence
for the pair $(A,L)$, in 1.1, ( since all the divisors in the linear system of
$M^2$, are symmetric for the action of $i$).
\erm

\begin{thebibliography}{99}

\bib [BL]{BL}  Birkenhake, Ch., Lange, H. : {\em Complex abelian varieties},   
 Springer-Verlag, Berlin, (1992).

\bib [I]{I} Iyer, J.:{\em Projective normality of abelian surfaces
given by primitive line bundles},  Manuscr.Math.$\bf{98}$, 139-153, (1999).

\bib [L]{L} Lazarsfeld, R.: {\em  Projectivite normale des surface  abeliennes}, Redige par O.Debarre.
Prepublication No. $\bf{14}$, Europroj- C.I.M.P.A., Nice, (1990).

\end {thebibliography}

\end{document}